\pdfoutput=1
\RequirePackage{ifpdf}
\ifpdf 
\documentclass[pdftex]{sigma}
\else
\documentclass{sigma}
\fi

\usepackage{tikz}
\usepackage{tikz-cd}

\numberwithin{equation}{section}

\newtheorem{Theorem}{Theorem}[section]
\newtheorem*{Theorem*}{Theorem}

\newtheorem{Lemma}[Theorem]{Lemma}

\theoremstyle{definition}
\newtheorem{Definition}[Theorem]{Definition}

\newtheorem{Remark}[Theorem]{Remark}

\begin{document}

\allowdisplaybreaks

\newcommand{\arXivNumber}{2507.13049}

\renewcommand{\PaperNumber}{078}

\FirstPageHeading

\ShortArticleName{A Note on Pliability and the Openness of the Multiexponential Map}

\ArticleName{A Note on Pliability and the Openness\\ of the Multiexponential Map in Carnot Groups}

\Author{Fr\'{e}d\'{e}ric JEAN~$^{\rm a}$, Mario SIGALOTTI~$^{\rm b}$ and Alessandro SOCIONOVO~$^{\rm a}$}

\AuthorNameForHeading{F.~Jean, M.~Sigalotti and A.~Socionovo}

\Address{$^{\rm a)}$~Unit\'{e} de Math\'{e}matiques Appliqu\'{e}es, ENSTA, Institut Polytechnique de Paris,\\
\hphantom{$^{\rm a)}$}~91120 Palaiseau, France}
\EmailD{\mail{frederic.jean@ensta.fr}, \mail{alessandro.socionovo@ensta.fr}}

\Address{$^{\rm b)}$~Laboratoire Jacques-Louis Lions, Sorbonne Universit\'{e}, Universit\'{e} de Paris,\\
\hphantom{$^{\rm b)}$}~CNRS, Inria, Paris, France}
\EmailD{\mail{mario.sigalotti@inria.fr}}

\ArticleDates{Received March 06, 2026, in final form August 11, 2026; Published online August 19, 2026}

\Abstract{In recent years, several notions of non-rigidity of horizontal vectors in Carnot groups have been proposed, motivated, in particular, by the characterization of monotone sets and Whitney extension properties. In this note we compare some of these notions.}


\Classification{53C17; 22E25; 93B05}

\section{Introduction}

Consider a Carnot group $\mathbb{G}$ with stratified Lie algebra ${\mathfrak g}={\mathfrak g}_1\oplus \dots \oplus {\mathfrak g}_s$ (see Definition~\ref{def:carnotgp}).
The endpoint map $E\colon L^\infty([0,1],{\mathfrak g}_1)\to \mathbb{G}$ (see Definition~\ref{def:endpoint})
 plays a crucial role in studying
metric and topological properties of $\mathbb{G}$. The goal of this paper is to explore and clarify the connections among several openness and submersivity properties of the endpoint map and of the multiexponential map
\begin{equation*} 
 \Gamma^{(p)}\colon \ ({\mathfrak g}_1)^p\to\mathbb{G}, \qquad \Gamma^{(p)}(Y_1,\dots,Y_p):=\exp(Y_{1})\cdots \exp(Y_{p}),
\end{equation*}
where $p\in \mathbb{N}$ and $\exp\colon {\mathfrak g}_1\to\mathbb{G}$ is the exponential map, i.e., the restriction of $E$ to ${\mathfrak g}_1$, seen as the subspace of $L^\infty([0,1],{\mathfrak g}_1)$ made of constant functions. In particular, given a vector $X\in{\mathfrak g}_1$, we are interested in the following properties of the multiexponential map and of $E_X(\cdot):=E(X+\cdot)$, the endpoint map based at $X$:
\begin{itemize}\itemsep=0pt\setlength{\leftskip}{0.22cm}
 \item[(H)] there exists $p\in\mathbb{N}$ such that the map $\Gamma^{(p)}$ is open at $(X,\dots,X)\in({\mathfrak g}_1)^p$;
 \item[(SbH)] there exists $p\in\mathbb{N}$ such that the map $\Gamma^{(p)}$ is a submersion at $(X,\dots,X)\in({\mathfrak g}_1)^p$;
 \item[(P)] the map $E$ is open at $X$ (i.e., $E_X$ is open at 0).
\end{itemize}

One of the earliest appearances of these conditions in the literature is the
observation by Cheeger and Kleiner in~\cite{CheegerKleiner} that in the three-dimensional Heisenberg group every $X\in {\mathfrak g}_1$ satisfies
condition (SbH)\ with $p=2$. As a consequence, they deduced the following \emph{cone property}: if a horizontally convex set $C$ contains a point $p$ in its closure and an open ball $B$, then it also contains
the cone over $B$ with vertex at $p$.

In the paper~\cite{ACM}, Arena, Caruso, and Monti generalized the submersivity of the map $(Y_1,Y_2)\mapsto \exp(Y_2) \exp(Y_1)$ to the case where $\mathbb{G}$ is of step two of M\'etivier type and
Morbidelli proved in~\cite{Morbidelli} that the cone property holds for all Carnot groups of step two.

Montanari and Morbidelli later observed
in~\cite{MontanariMorbidelli2021} that the cone property can be obtained
under the weaker assumption
that every $X\in {\mathfrak g}_1$ satisfies condition (H) (which is clearly implied by (SbH)). In addition, they showed that if a vector $X\in {\mathfrak g}_1$ satisfies condition (SbH), then the
Carnot--Carath\'eodory distance from the unit element of $\mathbb{G}$ is Pansu differentiable at $\exp(X)$.
(They also showed that (SbH)\ implies the related
notion of \emph{deformable direction} introduced by Pinamonti and Speight in~\cite{PinamontiSpeight}).

Condition (H) appears again in the paper~\cite{Rigot} by Rigot and it is used to characterize precisely monotone subsets of $\mathbb{G}$. We adapt the terminology introduced in~\cite{Rigot}, by distinguishing between the {\em {\rm (H)}-condition} and the {\em submersive {\rm (H)}-condition} (SbH)\ introduced above.

On the other hand, condition (P), called \emph{pliability}, was introduced
 in~\cite{Juillet-Sigalotti}, where Juillet and the second author proved the $C^1$ Whitney extension property for horizontal curves in Carnot groups where every vector in ${\mathfrak g}_1$ satisfies the property (P) (see also the recent preprint~\cite{SpeightZimmerman} by Speight and Zimmerman for a generalization to horizontal curves tangent to pliable directions).
The definition of pliability proposed in~\cite{Juillet-Sigalotti} is actually given in terms of the openness of the restriction of the endpoint map to a suitable proper subspace of $L^\infty([0,1],\mathfrak{g}_1)$. The equivalence between (P) and pliability as defined in~\cite{Juillet-Sigalotti} is one of the contributions of this paper. For a precise discussion on this point, see Remark~\ref{rem:plia-F}.

The first question addressed in this paper concerns the relationship
between pliability and the (H)-conditions.
Later, Sacchelli and the second author introduced in~\cite{SS18} the notion of {\em strong pliability {\rm (SP)}}, and they improved the results of~\cite{Juillet-Sigalotti} for horizontal curves in sub-Riemannian manifolds where every vector in ${\mathfrak g}_1$ of the tangent Carnot group satisfies (SP). Condition (SP) and the analogous notion of \emph{$p$-free {\rm (H)}-condition} (FH) are defined as follows:
\begin{itemize}\itemsep=0pt\setlength{\leftskip}{0.05cm}
 \item[(SP)] for all $\eta>0$ there exists $Y\in L^\infty([0,1],{\mathfrak g}_1)$ such that $\|Y\|_\infty<\eta$, $E_X(Y)=E_X(0)$, and the map $E_X$ is a submersion at $Y$;
 \item[(FH)] for all $\eta>0$ there exist an integer $p\ge 2$ and $W_1,\dots,W_p\in {\mathfrak g}_1$ such that $\|W_i-X\|<\eta$ for $i=1,\dots,p$, $\Gamma^{(p)} (W_1,\dots,W_p) = \Gamma^{(p)} (X, \dots, X)$, and $\Gamma^{(p)}$ is a submersion at $(W_1,\dots,W_p)$.
\end{itemize}

Finally, note that (SP) holds true in particular for $X\in{\mathfrak g}_1$ that are regular points of the~endpoint map
(i.e., those $X$ for which the straight line $t \mapsto \exp (tX)$ is not abnormal). This~property motivates a last condition:
\begin{itemize}\itemsep=0pt\setlength{\leftskip}{0.15cm}
 \item[(Reg)] the map $E$ is a submersion at $X$ (i.e., $E_X$ is a submersion at $0$).
\end{itemize}

Our main result is the following.

\begin{Theorem}\label{thm:main} The following hold:
\begin{itemize}\itemsep=0pt
 \item pliability, strong pliability,
 and the {\rm (FH)}-condition are equivalent;
 \item regularity {\rm (Reg)} and the submersive {\rm (H)}-condition are equivalent and imply the {\rm (H)}-con\-di\-tion, which, in turn, implies
 the preceding three conditions.
\end{itemize}
\end{Theorem}

We note that the equivalence between the (H)-condition and the submersive (H)-condition does not hold; see Remark~\ref{rem:nostrongH}.
Theorem~\ref{thm:main} is proved in Section~\ref{sec:proof} in a slightly more general version, see Theorem~\ref{thm:strong_main}.

\section{Preliminaries}

Let us recall some basic facts about Carnot groups, endpoint, and exponential maps.

\begin{Definition}\label{def:carnotgp}
 A Carnot group is a nilpotent, connected, and simply connected Lie group $\mathbb{G}$ whose Lie algebra ${\mathfrak g}$ admits a stratification, i.e., a direct sum decomposition{\samepage
 \begin{equation*}
 {\mathfrak g}={\mathfrak g}_1\oplus \dots \oplus {\mathfrak g}_s,
 \end{equation*}
 such that $[{\mathfrak g}_1,{\mathfrak g}_k]=
 {\mathfrak g}_{k+1}$ for $k=1,\dots,s-1$ and $[{\mathfrak g}_1,{\mathfrak g}_s]=0$.}
\end{Definition}

The \emph{left translation} by $g\in \mathbb{G}$ is the map $L_g\colon \mathbb{G}\to\mathbb{G}$ defined by $L_g(x):=gx$. Through this map, we have the identification
\begin{equation*}
 (L_g)_*\colon \ {\mathfrak g}\to T_g\mathbb{G}, \qquad {\mathfrak g}\ni X\mapsto (L_g)_*X\in T_g\mathbb{G}.
\end{equation*}
Every vector $X\in{\mathfrak g}$ defines the left-invariant vector field $X(g):=(L_g)_*X$. By abuse of notation, we identify the vector and the induced left-invariant vector field, in particular $X=X(1_\mathbb{G})$.

From now on, $I:=[0,1]\subset\mathbb{R}$ denotes the unit interval.
An absolutely continuous curve $\gamma\colon I\to\mathbb{G}$ is said to be \emph{horizontal} if
$\dot\gamma(t)\in (L_{\gamma(t)})_*{\mathfrak g}_1$
 for a.e.\ $t\in I$, that is,
\begin{equation*}
 \bigl(L_{\gamma(t)}^{-1}\bigr)_*
 \dot\gamma(t)=Y(t)\qquad\mbox{a.e.,}
\end{equation*}
for some function $Y\in L^\infty(I,{\mathfrak g}_1)$, called the \emph{control} of $\gamma$.
\begin{Definition}
\label{def:endpoint}
The \emph{endpoint map} is the mapping
 \begin{equation*}
 E\colon \ L^\infty(I,{\mathfrak g}_1)\to\mathbb{G}, \qquad E(Y):=\gamma(1),
 \end{equation*}
 where $\gamma$ is the unique solution to the Cauchy problem
 \begin{gather*}
 \bigl(L_{\gamma(t)}^{-1}\bigr)_* \dot\gamma(t)=Y(t),\qquad t\in I,\\
 \gamma(0)=1_\mathbb{G}.
 \end{gather*}
Given $X\in{\mathfrak g}_1$, the \emph{endpoint map based at $X$} is the mapping $E_X\colon L^\infty(I,{\mathfrak g}_1)\to \mathbb{G}$ defined by~${E_X(\cdot):=E(X+\cdot)}$.
\end{Definition}

In the following, we need the following well known result, which is proved, for instance, in~\cite[Section~8]{ABB20}.

\begin{Lemma}
 For every $X\in{\mathfrak g}_1$, the endpoint map $E_X\colon L^\infty(I,{\mathfrak g}_1)\to\mathbb{G}$ is of class $C^1$ as soon as $L^\infty(I,{\mathfrak g}_1)$ is endowed with the $L^2$ topology.
\end{Lemma}

If $X\in{\mathfrak g}$ and $\gamma\colon I \to\mathbb{G}$ is the integral curve of (the left-invariant vector field identified with)~$X$ with $\gamma(0)=1_\mathbb{G}$, then the {\em exponential map} is defined by
\begin{equation*}
 \exp\colon\ {\mathfrak g}\to\mathbb{G}, \qquad \exp(X):=\gamma(1).
\end{equation*}
Notice that $\gamma(t)=\exp(tX)$. Such a curve is called a {\em straight line}.

We denote the flow of a left-invariant vector field $X$ by ${\rm e}^{tX}\colon \mathbb{G}\to\mathbb{G}$. We have that
\[
{\rm e}^{tX}(g)=g\exp(tX).
\]

For every $\lambda\in\mathbb{R}$, we define the dilation on the Lie algebra $\delta_\lambda\colon {\mathfrak g}\to{\mathfrak g}$ as the linear map such~that
\[
 \delta_\lambda(X)=\lambda^i X \qquad \text{for all } X\in{\mathfrak g}_i.
\]
By abuse of notation, we use the same symbol to denote the dilation on the group $\delta_\lambda\colon \mathbb{G} \to \mathbb{G}$ defined by
\begin{equation}\label{eq:dilations-flow}
\delta_\lambda (\exp(X)) = \exp (\delta_\lambda(X))
\end{equation}
for every $X \in {\mathfrak g}$, so that $(\delta_\lambda)_* = \delta_\lambda$. In particular,
for
every $X \in {\mathfrak g}_1$ and $Y \in L^\infty(I,{\mathfrak g}_1)$, one has
\begin{equation}
\label{eq:dilations-flow2}
\delta_\lambda (E_X(Y)) = E_{\lambda X} (\lambda Y).
\end{equation}

\section{Proof of main result}\label{sec:proof}

Let us give a slightly more general form of the definitions of (P) and (SP).

\begin{Definition}\label{def:pl}
Let $V$ be a vector subspace of $L^\infty(I,{\mathfrak g}_1)$. A vector $X\in {\mathfrak g}_1$ is \emph{$V$-pliable} if the map $E_X$ restricted to $V$ is open at $0$.
$L^\infty(I,{\mathfrak g}_1)$-pliability coincides with pliability as defined by condition (P).
\end{Definition}

\begin{Definition}\label{def:strongpl}
 Let $V$ be a vector subspace of $L^\infty(I,{\mathfrak g}_1)$. The vector $X\in{\mathfrak g}_1$ is said to be \emph{strongly $V$-pliable} if for all $\eta>0$ there exists $Y\in V$ such that
 \begin{itemize}\itemsep=0pt
 	\item[(i)] $\|Y\|_\infty<\eta$;
 	\item[(ii)] $E_X(Y)=E_X(0)$;
 	\item[(iii)] $E_X|_V$ is a submersion at $Y$ (i.e., the differential of $E_X|_V$ at $Y$ is surjective).
 \end{itemize}
 $L^\infty(I,{\mathfrak g}_1)$-strong pliability coincides with strong pliability as defined by condition (SP).
\end{Definition}

\begin{Remark}\label{rem:plia-F}In~\cite[Definition 3.4]{Juillet-Sigalotti} and~\cite[Definition~4.1]{SS18}, the pliability
and the strong pliability
of
a vector $X$ are defined in terms of the openness
and the submersion property
of the
extended map $\mathcal{E}_X(Y):=(E_X(Y),Y(1))$ restricted to the space
\[
\mathcal C_0:=\bigl\{Y\in C^0(I,{\mathfrak g}_1)\mid Y(0)=0\bigr\}.
\]
However, in~\cite[Proposition~3.7]{Juillet-Sigalotti}
the openness of $\mathcal{E}_X|_{\mathcal C_0}$ is shown to be equivalent to the $\mathcal C_0$-pliability (in the sense of Definition~\ref{def:pl}), while in~\cite[Lemma~4.2]{SS18}
the strong pliability as defined in~\cite[Definition~4.1]{SS18} is shown to be
equivalent to that introduced in Definition~\ref{def:strongpl}.
As we will see in Theorem~\ref{thm:strong_main} (considering both $V={\mathcal C}_0$ and $V=L^\infty(I,{\mathfrak g}_1)$), it turns out that $X$ is $\mathcal C_0$-pliable if and only if it is pliable, since both properties are equivalent to the $p$-free (H)-condition.
\end{Remark}

Before stating Theorem~\ref{thm:main} in a slightly generalized version, let us introduce a useful notion.

\begin{Definition}\label{def:stabledensity}
We say that a subspace $V$ of $L^\infty(I,{\mathfrak g}_1)$ is
 \emph{$L^\infty$-stably $L^2$-dense in $L^\infty(I,{\mathfrak g}_1)$}
if there exists a constant $C>0$ such that, for every
$Y\in L^\infty(I,{\mathfrak g}_1)$, there exists a sequence
$(Z_k)_{k\in\mathbb{N}}$ in $V$ converging to $Y$ in $L^2(I,{\mathfrak g}_1)$ and such that
\[
\|Z_k\|_\infty
\le C\|Y\|_\infty
\qquad\text{for every }k\in\mathbb{N}.
\]
\end{Definition}
Examples of $L^\infty$-stably $L^2$-dense subspaces of $L^\infty(I,{\mathfrak g}_1)$
include $\mathcal{C}_0$ and the subspace of $L^\infty(I,{\mathfrak g}_1)$ made of piecewise constant functions.

Theorem~\ref{thm:main} is a corollary of the following result.

\begin{Theorem}\label{thm:strong_main}
 Let $\mathbb{G}$ be a Carnot group with Lie algebra ${\mathfrak g}$. Let $V$ be an
 $L^\infty$-stably $L^2$-dense subspace of $L^\infty(I,{\mathfrak g}_1)$.
 For every $X\in {\mathfrak g}_1$, the following conditions are equivalent:
 \begin{itemize}\itemsep=0pt\setlength{\leftskip}{0.5cm}
 \item[$V$-{\rm (P)}] $X$ is $V$-pliable;
 \item[$V$-{\rm (SP)}] $X$ is strongly $V$-pliable;
 \item[{\rm (FH)}] $X$ satisfies the $p$-free {\rm (H)}-condition.
 \end{itemize}
Moreover, $X$ satisfies the submersive {\rm (H)}-condition {\rm (SbH)} if and only if $X$ is a regular point of~$E$ {\rm ((Reg)}-condition$)$, and both conditions imply
\begin{itemize}\itemsep=0pt\setlength{\leftskip}{0.45cm}
\item[{\rm (H)}] $X$ satisfies the {\rm (H)}-condition;
\end{itemize}
which, in turn, implies the three equivalent conditions $V$-{\rm (P)}, $V$-{\rm (SP)}, and {\rm (FH)}.
\end{Theorem}

The structure of the proof of Theorem~\ref{thm:strong_main} is summarized
 in the following
 diagram:
\[
\begin{tikzcd}
 \text{(Reg)} \arrow[rr,Rightarrow] \arrow[d,dashed,Leftrightarrow]&&\mbox{(SP)}\arrow[d,Rightarrow]
 &&\mbox{(FH)}\arrow[ll,dashed,Rightarrow]
 \\
 \mbox{(SbH)}\arrow[r,Rightarrow]& \mbox{(H)} \arrow[r,dashed,Rightarrow]&\mbox{(P)}\arrow[r,Rightarrow]&V\mbox{-(P)}
 \arrow[r, dashed, Rightarrow]&V\mbox{-(SP)}
 \arrow[u,dashed,Rightarrow]
 \end{tikzcd}
\]
Here solid arrows denote implications that
 follow immediately from
the definitions, while dashed arrows are proved
via different lemmas.

The implication
\begin{tikzcd}[cramped,sep=small]V\mbox{-(P)}\arrow[r,dashed, Rightarrow]& V\mbox{-(SP)}\end{tikzcd}
 is a consequence of Lemmas \ref{lem:pl_iff_strpl} and \ref{lem:pl_iff_strpl2} (as explained in Remark~\ref{rmk:Vdense}). The implications \begin{tikzcd}[cramped,sep=small]\mbox{(H)}\arrow[r,dashed, Rightarrow]& \mbox{(P)}\end{tikzcd}
 and
 \begin{tikzcd}[cramped,sep=small]\mbox{(FH)}\arrow[r,dashed, Rightarrow]& \mbox{(SP)}\end{tikzcd}
 are proved in Lemma~\ref{lem:iii)-->i)}. The implication \begin{tikzcd}[cramped,sep=small]V\mbox{-(SP)}\arrow[r,dashed, Rightarrow]& \mbox{(FH)}\end{tikzcd} is proved in Lemma \ref{lem:strongly pliable implies weak-strong H}. The equivalence \begin{tikzcd}[cramped,sep=small]\text{(Reg)}\arrow[r,dashed, Leftrightarrow]& \text{(SbH)}\end{tikzcd} is proved in Lemma~\ref{le:strong=reg}. Finally, we notice that the implications (Reg)$\implies$(SP) and (SbH)$\implies$(FH) cannot be reversed, see Remark~\ref{rem:nostrongH}.

\begin{Lemma}\label{lem:pl_iff_strpl}
 A vector $X\in{\mathfrak g}_1$ is pliable if and only if it is strongly pliable.
\end{Lemma}

\begin{proof}
 The fact that strong pliability implies pliability is trivial. The converse can be proved in the following way.

Assume that $X\in{\mathfrak g}_1$ is pliable. It results from the homogeneity property~\eqref{eq:dilations-flow2} of $E_X$ that~\smash{$\frac X2$} and~\smash{$-\frac X2$} are also pliable. Let $g_1=\exp(X)$ and \smash{$g_{\frac12}=\exp\bigl(\frac X2\bigr)$}.
Fix $\eta>0$ and denote by~$B_\eta$ be the ball of radius $\eta$ in $L^\infty(I,{\mathfrak g}_1)$ centered at $0$.
By pliability of the vector~\smash{$-\frac X2$}, the set
${\mathcal U_-=g_1 \smash{E_{-\frac X2}(B_\eta)}}$ is a neighborhood of \smash{$g_{\frac12}$}.
 Now we can fix $0<\eta'<\eta$ such that~$\mathcal U:=\smash{E_{\frac X2}(B_{\eta'})\subset
 \mathcal U_-}$. Since \smash{$\frac X2$} is pliable, $\mathcal U$ is a neighborhood of \smash{$g_{\frac12}$}.

 Consider
 \begin{equation*}
 \bigl(L_{\gamma(s)}^{-1}\bigr)_*\dot\gamma(s)=\frac{X(\gamma(s))}{2}+Y_s(\gamma(s)),
\end{equation*}
 which can be seen as a control system
 with control $s\mapsto Y_s$ taking values in the ball of center $0$ and radius $\eta'$ in ${\mathfrak g}_1$.
 By a well-known result in geometric control theory (see, e.g.,~\cite[Section~3, Theorem~1]{Jurdjevic}),
 there exist $k\in \mathbb{N}$, $Y_1,\dots,Y_k\in {\mathfrak g}_1$ with $\|Y_i\|<\eta'$ for $i=1,\dots,k$, and ${t_1,\dots,t_k>0}$ such that
 \smash{\raisebox{-0.4pt}{$(s_1,\dots,s_k)\mapsto {\rm e}^{s_k(\frac X2+Y_k)}\circ \dots \circ {\rm e}^{s_1(\frac X2+Y_1)}(1_{\mathbb{G}})$}} is a submersion at $(t_1,\dots,t_k)$, with \smash{$\sum_{i=1}^kt_i$} arbitrary small. Since $t_k$ can be
 increased without affecting
 the submersion property, we can assume without loss of generality that \smash{$\sum_{i=1}^kt_i=1$}.

Set, for every $\varepsilon\in \mathbb{R}$ and $j\in\{1,\dots,k\}$,
\[\tilde Y_j(\varepsilon)=Y_j + \frac{\varepsilon}{t_j}\biggl( \frac X2+Y_j\biggr),\]
and notice that
\[t_j\biggl(\frac X2+\tilde Y_j(\varepsilon)\biggr)=(t_j+\varepsilon)\biggl(\frac X2+Y_j\biggr),\qquad j\in\{1,\dots,k\},\ \varepsilon\in \mathbb{R}.\]
As a consequence, the map
\[(\varepsilon_1,\dots,\varepsilon_k)\mapsto {\rm e}^{t_k(\frac X2+\tilde Y_k(\varepsilon_{k}))}\circ \dots \circ {\rm e}^{t_1(\frac X2+\tilde Y_1(\varepsilon_1))}(1_{\mathbb{G}})\]
is a submersion at $(\varepsilon_1,\dots,\varepsilon_k)=(0,\dots,0)$.
In particular, letting $Y'\in L^\infty(I,{\mathfrak g}_1)$
be equal to $Y_1$ on $[0,t_1)$, then to $Y_2$ on $[t_1,t_1+t_2)$, and so on, we have that $Y'\in B_{\eta'}$ and
\smash{$E_{\frac X2}$} is a submersion at $Y'$.

Letting \smash{$g'_{\frac12}=E_{\frac X2}(Y') \in \mathcal U$}, there exists $Y''\in B_\eta$ such that \smash{$g_1 E_{-\frac X2}(Y'')=g'_{\frac12}$}. We set
 \begin{equation*}
 Y(t)=
 \begin{cases}
 2Y'(2t), & t\in\bigl[0,\frac12\bigr],\\
 -2Y''(2(1-t)), & t\in\bigl[\frac12,1\bigr].
 \end{cases}
 \end{equation*}
 Then, we have $E_X(Y)=g_1$ and, since \smash{${\rm d}_{Y'}E_{\frac X2}$} is surjective, also the differential at~${\rm d}_YE_{X}$ is~surjective.
 Since $\|Y\|_\infty<2\eta$ and $\eta$ is arbitrarily small, the strong pliability of $X$ is proved.
\end{proof}

\begin{Remark}Let us stress that
Lemma~\ref{lem:pl_iff_strpl}
relies essentially on the special
properties of the endpoint map, in the sense that for the general case of a smooth open map $F$ at a point $X$, it is not true that there exist points arbitrarily close to $X$ whose image through $F$ coincide with $F(X)$ and at which $F$ is a submersion. Consider, for instance, \smash{$F\colon L^\infty(I,\mathbb{R})\ni X\mapsto \bigl(\int_I X(t){\rm d}t\bigr)^3\in \mathbb{R}$}. The map $F$ is open at $0$, but for every $X\in L^\infty(I,\mathbb{R})$ with $F(X)=0$ the differential of $F$ at $X$ is the null map.
\end{Remark}

\begin{Lemma} \label{lem:pl_iff_strpl2}
 Assume that
 $V$ is an
 $L^\infty$-stably $L^2$-dense subspace of
 $L^\infty(I,{\mathfrak g}_1)$.
 Then, a vector $X\in{\mathfrak g}_1$ is strongly
 pliable if and only if
 $X$ is strongly $V$-pliable.
\end{Lemma}

\begin{proof}
Since one implication is trivial, let us assume that $X$ is strongly
pliable and let us prove that it is strongly $V$-pliable.
 Fix $\eta>0$ and $C>0$ as in Definition~\ref{def:stabledensity}.
 Then there exists $Y\in L^\infty(I,{\mathfrak g}_1)$
 with $\|Y\|_\infty<\frac{\eta}{2C}$ such that $E_X(Y)=E_X(0)$ and
 ${\rm d}E_X(Y)$ is surjective. Pick $Y_1,\dots,Y_d\in L^\infty(I,{\mathfrak g}_1)$ such that
${\rm d}E_X(Y)Y_1,\dots,{\rm d}E_X(Y)Y_d$ is a basis of \smash{$T_{E_X(0)}\mathbb{G}$}, where $d$ is the dimension of $\mathbb{G}$.
 Since $V$ is dense in $L^\infty(I,{\mathfrak g}_1)$ for the $L^2$-topology,
 and the map
 $Z\mapsto {\rm d}E_X(Y)Z$ is continuous in $L^2(I,{\mathfrak g}_1)$, we can assume that $Y_1,\dots,Y_d\in V$
 by a standard approximation argument.

Consider the function $\mathbb{R}^d\ni (\alpha_1,\dots,\alpha_d)\mapsto E_X(Y+\alpha_1 Y_1+\dots+\alpha_d Y_d)$.
 By the inverse function theorem, for $r>0$ sufficiently small,
 such a function diffeomorphically maps $D=\{(\alpha_1,\dots,\alpha_d)\in \mathbb{R}^d \mid |\alpha_1|,\dots,|\alpha_d|<r \}$ into $E_X(D)$.
 Without loss of generality, \[r<\frac{\eta}{2(\|Y_1\|_\infty+\dots+\|Y_d\|_\infty)}.\]

 Let $Z\in V$ be such that $\|Z-Y\|_2<\varepsilon$, for some $\varepsilon>0$ to be fixed later.
 Since $V$ is
 $L^\infty$-stably $L^2$-dense in $L^\infty(I,{\mathfrak g}_1)$, we can assume that \smash{$\|Z\|_\infty\le C\|Y\|_\infty<\frac{\eta}2$}.
 Notice that
 ${E_X(Z+\alpha_1 Y_1+\dots+\alpha_d Y_d)}$ is arbitrarily close to $E_X(Y+\alpha_1 Y_1+\dots+\alpha_d Y_d)$ as $\varepsilon$ goes to zero, uniformly with respect to $(\alpha_1,\dots,\alpha_d)\in D$
 in any fixed system of coordinates.
 Using a classical degree argument similar to the one in~\cite[Lemma~5.5]{SS18}, and choosing $\varepsilon>0$ small enough, we have that $E_X(0)=E_X(w)$ for some $w\in \{Z+\alpha_1 Y_1+\dots+\alpha_d Y_d\mid |\alpha_1|,\dots,|\alpha_d|<r\}\subset V$ and ${\rm d}E_X(w)Y_1,\dots,{\rm d}E_X(w)Y_d$ is a basis of $T_{E_X(0)}\mathbb{G}$.

 Notice that, up to further reducing $\varepsilon$ if necessary,
 \[
 \|w\|_\infty\le \|Z\|_\infty +r (\|Y_1\|_\infty+\dots+\|Y_d\|_\infty)< \eta.
 \]
 Therefore, Properties (i), (ii), and (iii) in Definition \ref{def:strongpl} are satisfied by $w\in V$.
\end{proof}

\begin{Remark}\label{rmk:Vdense}
 Let $V$ be
 an $L^\infty$-stably $L^2$-dense subspace of $ L^\infty(I,{\mathfrak g}_1)$.
 Then, by the two previous results, it follows that
 \begin{equation*}
 \text{str. $V$-pliable} \implies \text{$V$-pliable} \implies \text{pliable} \stackrel{\rm Lemma \ \ref{lem:pl_iff_strpl}}{\implies} \text{str. pliable} \stackrel{\rm Lemma \ \ref{lem:pl_iff_strpl2}}{\implies} \text{str. $V$-pliable}.
 \end{equation*}
\end{Remark}

\begin{Lemma}\label{lem:iii)-->i)}
 If a vector $X\in{\mathfrak g}_1$ satisfies the {\rm (H)}-condition $($respectively, the {\rm (FH)}-condition$)$, then it is pliable $($respectively, strongly pliable$)$.
\end{Lemma}

\begin{proof}
If $X\in {\mathfrak g}_1$ satisfies the (H)-condition, there exists $p\in\mathbb{N}$ such that the map
 \begin{equation*}
 (Y_1,\dots,Y_p)\mapsto \exp(Y_{1}) \cdots \exp(Y_{p})
 \end{equation*}
 is open at $(X,\dots,X)$.
Hence, given any $\eta>0$,
$\bigl\{ {\rm e}^{Y_{p}}\circ \dots \circ {\rm e}^{Y_{1}}(1_{\mathbb{G}}) \mid \|Y_i-X\|<\eta\mbox{ for }i=1,\dots,p\bigr\}$
is a neighborhood of \smash{${\rm e}^{p X}(1_{\mathbb{G}})$}. Let \smash{$\delta_{\frac 1p}$} be the dilation in $\mathbb{G}$ of parameter $1/p$. Then, thanks to~\eqref{eq:dilations-flow}, for every $(Y_1,\dots,Y_p)\in ({\mathfrak g}_1)^p$ we have
\[{\rm e}^{{\frac 1p}Y_{p}}\circ \dots \circ {\rm e}^{{\frac 1p}Y_{1}}(1_\mathbb{G}) = \delta_{\frac 1p}\bigl({\rm e}^{Y_{p}}\circ \dots \circ {\rm e}^{Y_{1}}(1_\mathbb{G})\bigr).\]
Since, moreover, \smash{${\rm e}^{{\frac 1p}Y_p}\circ \dots \circ {\rm e}^{{\frac 1p}Y_1}(1_\mathbb{G})$} is the image by $E_X$ of the control constantly equal to $Y_i$ on $((i-1)/p,i/p)$ for $i=1,\dots,p$, it turns out that
the
image by $E_X$
of the set
\[\bigl\{Y\in L^{\infty}(I,{\mathfrak g}_1)\mid \|Y\|_\infty<\eta,\;Y|_{(i/p,(i+1)/p)}\mbox{ constant }\forall i=0,\dots,p-1\bigr\}\]
is
a neighborhood of ${\rm e}^{X}(1_\mathbb{G})$, which implies that $X$ is pliable.

The proof in the case of the (FH)-condition is similar.
\end{proof}

\begin{Lemma}\label{lem:strongly pliable implies weak-strong H}
If $X\in {\mathfrak g}_1$ is strongly pliable, then it satisfies the {\rm (FH)}-condition.
\end{Lemma}

\begin{proof}
Let $X\in{\mathfrak g}_1$ be strongly pliable, and denote by ${\rm PC}={\rm PC}(I,{\mathfrak g}_1)$ the subspace of $L^\infty(I,{\mathfrak g}_1)$ made of piecewise constant functions with rational discontinuity points. By Remark~\ref{rmk:Vdense}, $X$ is strongly ${\rm PC}$-pliable.

 Fix $\eta>0$. By the strong ${\rm PC}$-pliability, there exist $Z,Y_1,\dots,Y_d\in {\rm PC}$, with $\|Z\|_\infty<\eta$, such that $E_X(Z) = E_X(0)$ and ${\rm d}E_X(Z)Y_1,\dots,{\rm d}E_X(Z)Y_d$ is a basis of $T_{E_X(0)}\mathbb{G}$. By definition of the space ${\rm PC}$, there exists an integer $p\ge 2$ such that \smash{$Z|_{((i-1)/p,i/p)}$}, \smash{$Y_1|_{((i-1)/p,i/p)},\dots,Y_{d}|_{((i-1)/p,i/p)}$} are constant for every $i=1,\dots,p$.

 Denote by $Z_i$ and $Y_{j,i}$ the values of $Z$ and $Y_j$ on the interval \smash{$\bigl(\frac{i-1}{p}, \frac ip\big)$}. Then the map $\Gamma^{(p)}\colon ({\mathfrak g}_1)^p\to\mathbb{G}$ defined by
 \[
 \Gamma^{(p)}(W_1,\dots,W_p):= {\rm e}^{{W_{p}}}\circ \dots \circ {\rm e}^{{W_{1}}}(1_\mathbb{G}),
 \]
satisfies $\Gamma^{(p)} (X+Z_1,\dots,X+Z_p) = \Gamma^{(p)} (X, \dots, X)$.
 It remains to show
that $\Gamma^{(p)}$ is a submersion at $(X+Z_1,\dots,X+Z_p)$.

We have that \smash{$Y_j=\sum_{i=1}^p Y_{j,i}\chi _{((i-1)/p,i/p)}$} and so
 \[{\rm d}E_X(Z)Y_j=\sum_{i=1}^p {\rm d}E_X(Z)Y_{j,i}\chi _{(\frac{i-1}p,\frac ip)}.\]
Moreover,
\begin{gather}
 {\rm d}E_X(Z)Y_{j,i}\chi _{(\frac{i-1}p,\frac ip)}\nonumber\\
 \qquad{}=\frac{{\rm d}}{{\rm d}\varepsilon}\biggr|_{\varepsilon=0}
 \Bigl({\rm e}^{\frac{X+Z_p}p}\circ\dots\circ {\rm e}^{\frac{X+Z_{i+1}}p}\circ {\rm e}^{\frac{X+Z_i+\varepsilon Y_{j,i}}p}\circ {\rm e}^{\frac{X+Z_{i-1}}p}\circ \dots\circ {\rm e}^{\frac{X+Z_1}p}(1_{\mathbb{G}})\Bigr)\nonumber\\
 \qquad{}\stackrel{\eqref{eq:dilations-flow}}{=}\frac{{\rm d}}{{\rm d}\varepsilon}\biggr|_{\varepsilon=0}
 \delta_{\frac1p}\bigl({\rm e}^{X+Z_p}\circ\dots\circ {\rm e}^{X+Z_{i+1}}\circ {\rm e}^{X+Z_i+\varepsilon Y_{j,i}}\circ {\rm e}^{X+Z_{i-1}}\circ \dots\circ {\rm e}^{X+Z_1}(1_{\mathbb{G}})\bigr)\nonumber\\
\qquad{}=\bigl(\delta_{\frac1p}\bigr)_*\frac{{\rm d}}{{\rm d}\varepsilon}\biggr|_{\varepsilon=0}\bigl({\rm e}^{X+Z_p}\circ\dots\circ {\rm e}^{X+Z_{i+1}}\circ {\rm e}^{X+Z_i+\varepsilon Y_{j,i}}\circ {\rm e}^{X+Z_{i-1}}\circ \dots\circ {\rm e}^{X+Z_1}(1_{\mathbb{G}})\bigr)\nonumber\\
 \qquad{}=\bigl(\delta_{\frac1p}\bigr)_*\frac{{\rm d}}{{\rm d}\varepsilon}\biggr|_{\varepsilon=0}
 \Gamma^{(p)}(X+Z_{1},\dots,X+Z_{i-1},X+Z_i+\varepsilon Y_{j,i},X+Z_{i+1}, \dots, X+Z_{p})\nonumber\\
\qquad{}=\bigl(\delta_{\frac1p}\bigr)_* {\rm d}\Gamma^{(p)}(X+Z_1,\dots,X+Z_p)\cdot (0,\dots,0, Y_{j,i},0,\dots,0).\label{eq:diffdilated}
\end{gather}
So
the image of the differential of
$\Gamma^{(p)}$
at $(X\!+\!Z_1,\dots,X\!+\!Z_p)$ contains the image by \smash{$(\delta_{p})_*\!=\!\bigl(\delta_{\frac1p}\bigr)_*^{-1}$} of
the
span of ${\rm d}E_X(Z)Y_1,\dots,{\rm d}E_X(Z)Y_d$.
Hence $\Gamma^{(p)}$
is a submersion at $(X\!+\!Z_1,\dots,X\!+\!Z_p)$.
\end{proof}

\begin{Lemma}
\label{le:strong=reg}
A vector $X\in{\mathfrak g}_1$ satisfies the submersive {\rm (H)}-condition {\rm (SbH)}
 if and only if it is a regular point of the endpoint map~$E$.
\end{Lemma}
\begin{proof}
On the one hand, if $X\in{\mathfrak g}_1$ satisfies
(SbH), then a computation analogous to \eqref{eq:diffdilated} with $Z=0$ shows that ${\rm d}E_X(0)$ is onto.
On the other hand,
note that the argument of Lemma~\ref{lem:pl_iff_strpl2} shows that, if $E_X$ is a submersion at $0$, then $E_X$ restricted to PC is also a submersion at $0$. Then the argument used in the proof of Lemma~\ref{lem:strongly pliable implies weak-strong H} shows that $X$ satisfies the submersive (H)-condition.
\end{proof}

\begin{Remark}
 \label{rem:nostrongH}
The (H)-condition is strictly weaker than the submersive (H)-condition (SbH).
In particular, also pliability, strong pliability, and the (FH)-condition are strictly weaker than (SbH).
Indeed,
Carnot groups of step two
are known to satisfy the (H)-condition~\cite[Theorem~2.1]{Morbidelli},
 while non-M\'{e}tivier step two Carnot groups admit abnormal curves of the form $t\mapsto \exp(tX)$
 (see~\cite[Proposition~3.6]{Montanari-Morbidelli-JMMA})
 and, hence, vectors $X$ that do not satisfy the submersive (H)-condition.
\end{Remark}

\subsection*{Acknowledgements}
 We are thankful to
S\'everine Rigot for
 several
useful discussions
on
the topics discussed in this paper, and, in particular, for pointing out a mistake that we made in a previous version of this work.
We also thank the anonymous referees for their valuable remarks that greatly improved the paper.
This research benefited from the support of the FMJH.
 This project has received funding from the
{\em Fondazione Ing. Aldo Gini} of Padova, Italy, and from the
European Union’s Horizon 2020 research and innovation programme under the Marie Sk{\l}odowska-Curie grant agreement No 101034255. \includegraphics[width=0.65cm]{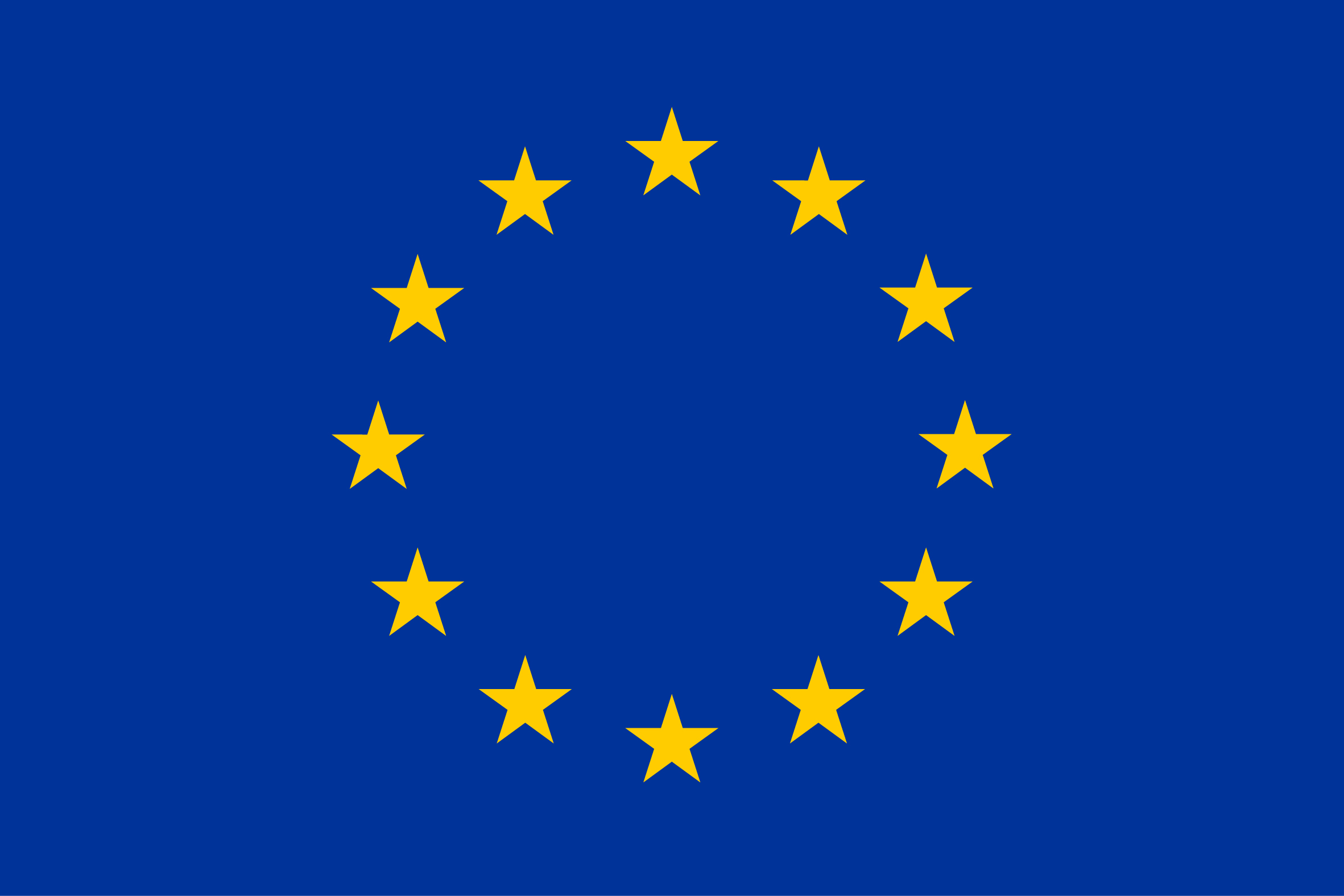}


\pdfbookmark[1]{References}{ref}
\LastPageEnding

\end{document}